\documentclass[11pt,twoside,a4paper]{article}
\usepackage{amsfonts}
\usepackage{amssymb}
\usepackage{amsmath}
\usepackage[french, english]{babel}

\begin{document}

\newcommand{\pa}{\partial}
\newcommand{\opa}{\overline\pa}
\newcommand{\ol}{\overline }

\numberwithin{equation}{section}

\newcommand\C{\mathbb{C}}  
\newcommand\R{\mathbb{R}}
\newcommand\Z{\mathbb{Z}}
\newcommand\N{\mathbb{N}}
\newcommand\PP{\mathbb{P}}

{\LARGE \centerline{Lorentzian $CR$ structures and nonembeddability}}
\vspace{0.8cm}

\centerline{\textsc {Judith Brinkschulte\footnote{Universit\"at Leipzig, Mathematisches Institut, Augustusplatz 10, D-04109 Leipzig, Germany. 
E-mail: brinkschulte@math.uni-leipzig.de}}
 and C. Denson Hill
\footnote{Department of Mathematics, Stony Brook University, Stony Brook NY 11794, USA. E-mail: dhill@math.stonybrook.edu\\
{\bf{Key words:} }abstract $CR$ structures, local embeddability, Lorentzian $CR$ manifolds\\
{\bf{2010 Mathematics Subject Classification:}}32V05, 32V30, 32G07}}

\vspace{0.5cm}

\begin{abstract} In this paper we construct examples of $CR$ deformations of Lorentzian hypersurfaces which are $CR$ embeddable at all points outside
an arbitrarily small compact set whose interior contains a point where $CR$
embeddablity is not possible.
\end{abstract}

\section{Introduction}

Let $M$ be a smooth differentiable manifold of odd dimension 
$\mathrm{dim}_{\R}M = 2n+1$. An abstract $CR$ structure of hypersurface type on $M$ is given
 by the choice of a complex subbundle $T^{0,1}M \subset\!\C\!
 \otimes TM$ of rank $n$ satisfying
$$ T^{1,0}M \cap T^{0,1}M = \lbrace 0 \rbrace \qquad 
\mathrm{where} \ T^{1,0}M = \ol{T^{0,1}M}$$
together with the integrability condition
$$ \lbrack\, \mathcal{C}^{\infty} (M, T^{0,1}M) ,  
\mathcal{C}^{\infty} (M, T^{0,1}M)\,\rbrack \subset  \mathcal{C}^
{\infty} (M, T^{0,1}M)$$
For every smooth real 
hypersurface $M$ in
$\C^{n+1}$,  the $n$-dimensional subbundle $T^{0,1}M
= T^{0,1}\C^{n+1} \cap \C\! \otimes\! TM$ defines the standard $CR$ 
structure on $M$.
A fundamental problem in the theory of $CR$ manifolds is to decide when an
abstract
$CR$ structure can be $CR$ embedded into some complex manifold. \\

In \cite{BH} we considered compactly supported $CR$ deformations of the standard $CR$ structure on quadratic $CR$ submanifolds of arbitrary $CR$ codimension. We showed that any such compactly supported $CR$ deformation is $CR$ embeddable if the quadratic $CR$ submanifold is $2$-pseudoconcave. 
But what if the quadratic $CR$ manifold is only
$1$-pseudoconcave? On the one hand, our method of proof, and the estimates
we used, do not appear to extend to that situation, even for codimension
one. On the other hand, we have not managed to construct a compactly
supported global formally integrable deformation which is not $CR$ embeddable in order to produce a counterexample. The difficulty in
doing so is clearly related to the problem encountered in \cite{JT1};
see \cite{JT2} and $\cite[\mathrm{Remark}\ 4.4]{HN2}$. In \cite{HN2} it was speculated that
in this "Lorentzian" situation, there might be geometric obstructions to
the existence of perturbations with arbitrarily small compact support
which are not $CR$ embeddable, coming from some propagation phenomena along
the "light rays". In this paper we construct examples of $CR$ deformations of Lorentzian hypersurfaces which are $CR$ embeddable at all points outside
an arbitrarily small compact set whose interior contains a point where $CR$
embeddablity is not possible. Thus a non embeddable point is "trapped"
and surrounded by embeddable points. This shows that a naive approach to
any such propagation phenomena does not appear to work.\\

More precisely, in this paper we consider the following set-up:

Let $N\subset \C^2$ be a strictly pseudoconvex hypersurface. We may assume that $0\in N$ and that for some small ball $V$ centered at $0$ in $\C^2$, $N\cap V$ is given by $N\cap V= \lbrace \zeta\in V\mid r(\zeta) =0\rbrace$ in a neighborhood of $0$, where $r$ is a strictly convex function in a neighborhood of $0$ satisfying $dr(0)\not= 0$\textsl{}. We consider a sufficiently small open neighborhood $U$ of $0\in\C^{n-1}$ and define
$$M = \lbrace (z,\zeta) \in U\times V \mid \varrho(z,\zeta)=0\rbrace,$$
where $\varrho(z,\zeta)= \vert z_1\vert^2 + \ldots +\vert z_{n-1}\vert^2 - r(\zeta)$.
Then, if $U$ and $V$ are sufficiently small, $M$ is a smooth Lorentzian hypersurface in $\C^{n+1}$ with its standard $CR$ structure: its Levi form has one negative and $n-1$ positive eigenvalues.
Let  $K\subset M$ be an arbitrarily small compact neighborhood of the origin. Then we can prove the following:

\newtheorem{main}{Theorem}[section]
\begin{main}   \label{main}   \ \\
There exists a family of abstract Lorentzian $CR$ structures $(T^{0,1}M_a)_{a>0}$ on $M$  with the following properties:
\begin{enumerate}
\item The abstract $CR$ structures $(T^{0,1}M_a)_{a>0}$ converge to the standard $CR$ structure on $M$ as $a$ tends to $0$ in the usual $\mathcal{C}^\infty$ topology.
\item The abstract $CR$ structures $(T^{0,1}M_a)_{a>0}$ coincide with the standard $CR$ structure on $M$ to infinite order at $0$ for all $a > 0$.
\item The abstract $CR$ structure $T^{0,1}M_a$ is not locally $CR$ embeddable at $0$ for $a > 0$.
\item There exists a $CR$ embedding of the abstract $CR$ structure $T^{0,1}M_a$ on $M\setminus K$ into some complex manifold of dimension $n+
1$.
\end{enumerate}
\end{main}

\section{Proof of the theorem}

We assume that the euclidean ball $V$ is small enough such that $N$ divides $V$ into exactly two open subsets and that $N$ intersects $\pa V$ transversally. We consider a small perturbation $\tilde r$ of $r$ which is still strictly plurisubharmonic on $V$ with the following properties: $\tilde r = r$ in some open neighborhood of $0$ and $\tilde r > r$ outside some ball centered at $0$. We obtain a small perturbation of $N$: $\tilde{N} = \lbrace \zeta\in V\mid \tilde{r}(\zeta)=0\rbrace$.\\

 We then define the two  open subsets $\Omega^+ = \lbrace \zeta\in  V\mid \tilde r(\zeta) > 0\rbrace$ and
$\Omega^- = \lbrace \zeta\in  V\mid \tilde r(\zeta) < 0\rbrace$.
In order to define the abstract $CR$ structures $T^{0,1}M_a$, we need a $\opa$-closed $(0,1)$-form $\omega$ on $\Omega^+$ which is not $\opa$-exact near $0$.\\

It follows from \cite{AH} that there exists a $\opa_N$-closed $(0,1)$-form $f$ defined in a neighborhood of $0$ on $N$ which is not $\opa_N$-exact on any open neighborhood of $0$ in $N$. Multiplying $f$ by a cutoff function with support in an arbitrary, but fixed compact $S$ in $N$, we may assume that $f$ is compactly supported in $N\cap V$, and still $f$ is not $\opa_N$-exact on any open neighborhood of $0$ in $N$. By choosing $S$ small enough, we may assume also assume that $\lbrace (z,\zeta)\in M\mid \zeta \in S\rbrace\subset K$.\\

It is well known that we have a long exact sequence induced from the Mayer-Vietoris sequence
$$ \ldots\longrightarrow H^{0,1}(\ol V)\longrightarrow H^{0,1}(\ol{\Omega^+})\oplus H^{0,1}(\ol{\Omega^-})\longrightarrow H^{0,1}(\tilde{N})\longrightarrow H^{0,2}(\ol V)\longrightarrow \ldots ....$$
The first and last cohomology group vanish since $V$ is a euclidean ball. This implies that 
$$f = (F^+ - F^-)_{\mid \tilde N},$$
with $F^\pm\in \mathcal{C}^\infty_{0,1}(\ol{\Omega^\pm})\cap\mathrm{Ker}\opa$. Since $\ol{\Omega^-}$ is the transversal intersection of strictly pseudoconvex domains, there exists $U\in\mathcal{C}^\infty_{0,0}(\ol{\Omega^-})$ satisfying $\opa U = F^-$ on $\ol{\Omega^-}$. But then $F^+_{\mid \tilde N} = f + \opa_{\tilde N} U_{\mid \tilde N}$, which implies that the restriction of $F^+$ to $N$ is not $\opa_N$-exact on any open neighborhood of $0$ in $N$.\\

By intersection of $\ol{\Omega^+}$ with the tangent plane to $\tilde N$ at $0$, we obtain a domain $D$ which is the transversal intersection of weakly pseudoconvex domains such that $\ol D \cap \tilde N = \lbrace 0\rbrace$. By \cite{D}, there exists a function $g\in \mathcal{C}^\infty_{0,0}(\ol D)$ satisfying $\opa g = F^+$ on $\ol D$. Let $\tilde g$ be a smooth extension of $g$ to $\ol{\Omega^+}$, and set $\omega = F^+ - \opa\tilde g$. Then $\omega= \omega_1 d\ol \zeta_1 + \omega_2 d\ol \zeta_2 $ is a $\opa$-closed $(0,1)$-form $\omega$ on $\Omega^+$ which is not $\opa$-exact near $0$. Moreover, $\omega$ vanishes to infinite order at $0$.\\

For $a> 0$ we consider the system 

\begin{equation}  \label{system}
\left\{
\begin{aligned}
\ol L_1 \quad & = & \frac{\pa}{\pa\ol z_1}, \hspace{2.5cm}\\
& \cdots & \\
\ol L_{n-1} & = & \frac{\pa}{\pa\ol z_{n-1}}, \hspace{2.3cm}\\
\ol L_n \quad & = & \frac{\pa}{\pa\ol \zeta_1} + a \omega_1(\zeta) \frac{\pa}{\pa z_{n-1}},\\
\ol L_{n+1} & = & \frac{\pa}{\pa\ol \zeta_2} + a \omega_2(\zeta) \frac{\pa}{\pa z_{n-1}}
\end{aligned}  \right.
\end{equation}

Clearly $\lbrack \ol L_j, \ol L_k\rbrack =0$ for all $k$ if $1\leq j\leq n-1$. And
$$\lbrack \ol L_n,\ol L_{n+1}\rbrack = a\lbrace \frac{\pa\omega_2}{\pa\ol\zeta_1} -\frac{\pa\omega_1}{\pa\ol\zeta_2}\rbrace \frac{\pa}{\pa z_{n-1}} = 0$$
since $\opa\omega =0$. Thus (\ref{system}) defines an integrable almost complex structure on $U\times \Omega^+\subset\C^{n+1}$, which is smooth up to the partial boundary $U\times\lbrace\tilde r = 0\rbrace$. By the Newlander-Nirenberg theorem, we thus have a new complex structure on $U\times \Omega^+$. The abstract $CR$ structure $T^{0,1}M_a$ on $M$ is defined as the induced $CR$ structure on $M$. For $a>0$ sufficiently small, its Levi form has again Lorentzian signature.\\

1. is clear in view of (\ref{system}). \\

Since $\omega$ vanishes to infinite order at $0$ we get 2.\\

 To justify 4. we note that at points on $M$ outside $K$, the $CR$ structure comes in fact from a complex structure defined in a neighborhood of that point.\\

The proof of 3. is as in \cite{HN}.\\

Indeed, on $M$ we may use the real coordinates $t,\xi$ with $t = (t_1, t_2, \ldots, t_{2n-2})$, $\xi = (\xi_1,\xi_2,\xi_3)$, where $z_1 = t_1 + t_2$, $z_2 = t_3 + i t_4$, $\ldots$, $ z_{n-1} = t_{2n-3}+ i t_{2n-2}$, $\zeta_1 = \xi_1 + i\xi_2$, $\zeta_2 = \xi_3 + i\xi_4$, are restricted to $M$.\\

Suppose that we have a local $CR$ embedding of $T^{0,1}M_a$ near the origin by $CR$ functions $u_1(t,\xi), u_2(t,\xi),\ldots, u_{n+1}(t,\xi)$ with $du_1\wedge\ldots\wedge du_{n+1}\not=0 $ at $(0,0)$. 

Each $u_j$ satisfies $\ol L_k u_j =\frac{\pa u_j}{\pa \ol z_k}=0$, for $k=1,\ldots, n-1$. The Riemann removability theorem and the H. Lewy two-sided extension of $CR$ functions then imply that $u_j$ is in fact holomorphic in $z_k$ for $k=1,\ldots, n-1$.
 Therefore the Jacobian of the embedding map has a block decomposition of the form

\begin{equation} \nonumber
\left\lbrack
\begin{aligned}
 \frac{\pa u}{\pa z} & \quad 0 & \frac{\pa u}{\pa \xi} \\
  0 \quad  & \frac{\pa \ol u}{\pa \ol z} & \frac{\pa \ol u}{\pa \xi}
\end{aligned}
\right\rbrack,
\end{equation}
where we write $u = (u_1, \ldots, u_{n+1})^T$ and $z = (z_1,\ldots, z_{n-1})$ for short. This matrix has $2(n+1)$ rows and $2n+1$ linearly independent columns. Hence $\frac{\pa u_j}{\pa z_{n-1}}\not= 0$ for some $j$ at $(0,0)$. By renaming, we may assume $\frac{\pa u_{n-1}}{\pa z_{n-1}}\not= 0$.\\

The coordinates on $\mathbb{C}^{n+1}$ also define $CR$ functions on $T^{0,1}M_a$:\\
 $z_1(t,\xi),\ldots, z_{n-1}(t,\xi), \zeta_1(t,\xi),\zeta_2(t,\xi)$, and $ dz_1\wedge\ldots\wedge dz_{n-2}\wedge du_{n-1} \wedge d\zeta_1\wedge d\zeta_2\not= 0$ at $(0,0)$. So we arrive at a new local embedding map
$$\varphi: (t,\xi) \mapsto ( z_1(t,\xi),\ldots, z_{n-2}(t,\xi), u_{n-1}(t,\xi), \zeta_1(t,\xi),\zeta_2(t,\xi))$$
of some neighborhood $W$ of $(0,0)$ into $\mathbb{C}^{n+1}$. $\varphi(W)$ is a piece of a real hypersurface  in $\mathbb{C}^{n+1}$. Let $w =(w_1,\ldots, w_{n+1})$ denote the coordinates in $\mathbb{C}^{n+1}$, and consider, for points on $\varphi(W)$, the function
$$F(w) = \varphi_{\ast}\big(- \lbrack \frac{\pa u_{n-1}}{\pa z_{n-1}}\rbrack^{-1}\big),$$
where $\varphi_{\ast}$ is the push-forward by the diffeomorphism $\varphi$ of $W$ onto $\varphi(W)$. It follows that $F$ is a $CR$ function on $\varphi(W)$; in particular, it is holomorphic in $w_{n-1}$ by the inverse mapping theorem for holomorphic functions of one variable. On $\varphi(W)$ we may define the function
\begin{equation}  \label{function}
G(w)= \int_0^{w_{n-1}} F(w_1,\ldots, w_{n-2}, \eta, w_n, w_{n+1}) d\eta
\end{equation}
by a contour integral in the $w_{n-1}$-plane. This is well defined by the open mapping theorem from one complex variable. We now pull back to get a function $g(t,\xi)= \varphi^\ast G$ on $V$, which is a $CR$ function there. This can be seen by replacing $F$ in (\ref{function}) by a smooth extension $\tilde F$ of $F$ off of $\varphi(W)$ such that $\opa\tilde F_{\mid\varphi(V)} =0$ and differentiating. Next we have
\begin{equation}  \nonumber
\begin{aligned}
\frac{\pa g}{\pa z_{n-1}} & = & \frac{\pa G}{\pa z_{n-1}}( w_1,\ldots, w_{n-2}, u_{n-1},w_n,w_{n+1}) \frac{\pa u_{n-1}}{\pa z_{n-1}}(z,\xi)\hspace{3cm}\\
\quad & = & F( w_1,\ldots,w_{n-2}, u_{n-1}(z,\xi),w_n, w_{n+1}) \frac{\pa u_{n-1}}{\pa z_{n-1}}(z,\xi)\hspace{3cm}\\
\quad & = & -1,\hspace{10cm}
\end{aligned}
\end{equation}
so $g(z) = -z_{n+1} + \chi(z_1,\ldots, z_{n-2}, \xi)$, where $\chi(z_1,\ldots, z_{n-2}, \xi)$ is a smooth "constant of integration". Now the fact that $g$ is a $CR$ function implies that $\ol L_k g =0$, hence $\ol L_k \chi- a\omega_k =0$ for $k=n,n+1$. But for $a\not= 0$, this means that there exists a neighborhood $V^\prime$ of $0$ on $N$ such that $\omega$ is $\opa_N$-exact on $V^\prime$. This is a contradiction to the assumption on $\omega$. Therefore for $a\not=0$, $T^{0,1}M_a$ is not locally $CR$ embeddable over any open neighborhood of $0$ on $M$. \hfill$\square$\\

\vspace{1cm}

{\bf Acknowledements.} The first author was supported by Deutsche Forschungsgemeinschaft (DFG, German Research Foundation, grant BR 3363/2-1).\\


\begin{thebibliography}
\footnotesize

\bibitem[AH]{AH} \textsc{A. Andreotti, C.D. Hill:} \emph{E.E. Levi convexity and the Hans Lewy problem. Part II: Vanishing theorems.} Ann. Sc. Norm. Sup. Pisa {\bf 26}, 747--806 (1972).

\bibitem[BH]{BH} \textsc{J. Brinkschulte, C. D. Hill:} \emph{Inflexible $CR$ submanifolds}, Math. Z. (2016),
doi:10.1007/s00209-016-1831-6

\bibitem[D]{D} \textsc{A. Dufresnoy:} \emph{Sur l'op\'erateur $d^{\prime\prime}$ et les fonctions diff\'erentiables au sens de Whitney.} Ann. Inst. Fourier {\bf 29}, 229--238 (1979).

\bibitem[H]{H} \textsc{C.D. Hill:} \emph{Counterexamples to Newlander-Nirenberg up to the boundary.} Proc. Symp. Pure Math. {\bf 52}, 191--197 (1991).

\bibitem[HN1]{HN} \textsc{C.D. Hill, M. Nacinovich:} \emph{Embeddable $CR$ manifolds with nonembeddable smooth boundary.} B.U.M.I. {\bf 7}, 387--395 (1993).

\bibitem[HN2]{HN2} \textsc{C.D. Hill, M. Nacinovich:} \emph{Non completely solvable systems of complex first order PDE's.} Rend. Sem. Mat. Univ. Padova {\bf 129}, 129--168 (2013).

\bibitem[JT1]{JT1} \textsc{H. Jacobowitz, F. Tr\`{e}ves:} \emph{Abberant $CR$ structures.} Hokkaido Math. J. {\bf 12}, 276--292 (1983).

\bibitem[JT2]{JT2} \textsc{H. Jacobowitz, F. Tr\`{e}ves:} \emph{Erratum: Abberant $CR$ structures.} Hokkaido Math. J. {\bf 42}, 473--474 (2013).


\end{thebibliography}
\end{document}